\documentclass[10pt]{article}
\usepackage{mathrsfs}
\usepackage{amsmath}
\usepackage{amssymb}
\usepackage{graphicx}
\def \dis{\hbox{dis}}
\def \la{\lambda}

\newtheorem{theorem}{\scshape \mdseries  Theorem}[section]
\newtheorem{lemma}[theorem]{\scshape \mdseries  Lemma}

\begin{document}

\title{\sf Minimum Distance Spectral Radius of Graphs with Given Edge Connectivity\thanks{
      Supported by National Natural Science Foundation of China (11071002), Program for New Century Excellent
Talents in University, Key Project of Chinese Ministry of Education
(210091), Specialized Research Fund for the Doctoral Program of
Higher Education (20103401110002), Science and Technological Fund of
Anhui Province for Outstanding Youth (10040606Y33), Scientific
Research Fund for Fostering Distinguished Young Scholars of Anhui
University(KJJQ1001), Academic Innovation Team of Anhui University
Project (KJTD001B), Fund for Youth Scientific Research of Anhui
University(KJQN1003). }}
\author{Xiao-Xin Li$^{1,2}$, Yi-Zheng Fan$^{2,}$\thanks{Corresponding author. E-mail addresses:
fanyz@ahu.edu.cn (Y.-Z. Fan), lxx@czu.edu.com (X.-X. Li),
wangy@ahu.edu.cn (Y. Wang)}, Yi Wang$^2$
\\
  {\small  \it $^1$Department of Mathematics, Chizhou University, Chizhou 247000, P.R. China} \\
  {\small \it $^2$School of Mathematical Sciences, Anhui University, Hefei 230601, P. R. China}
 }
\date{}
\maketitle

\noindent {\bf Abstract}\ \ In this paper we determine the unique graph with minimum distance spectral radius
among all connected graphs of fixed order and given edge connectivity.

\noindent {\bf MR Subject Classifications:} 05C50, 15A18


\section{Introduction}
Let $G$ be a connected simple graph with vertex set $V(G)$ and edge set $E(G)$.
The {\it distance} between two vertices $u,v $ of $G$, denoted by $\dis_{uv}$, is defined as the length of the shortest path between $u$ and $v$ in $G$.
The {\it distance matrix} of $G$, denoted by $D(G)$, is defined by $D(G)=(\dis_{uv})_{u,v \in V(G)}$.
 Since $D(G)$ is symmetric, its eigenvalues are all real. In addition, as $D(G)$ is nonnegative and irreducible, by Perron-Frobenius theorem,
 the spectral radius $\rho(G)$ of $D(G)$ (called the {\it distance spectral radius} of $G$), is exactly the largest eigenvalue of $D(G)$ with multiplicity one;
  and there exists a unique (up to a multiple) positive eigenvector corresponding to this eigenvalue, usually referred to the {\it Perron vector} of $D(G)$.

The distance matrix is very useful in different fields, including the design of communication networks \cite{gra1}, graph embedding theory \cite{ede,gra3,gra2} as well as molecular stability \cite{hos,rou}.
Balaban et al. \cite{bal} proposed the use of the distance spectral radius as a molecular descriptor.
 Gutman et al. \cite{gut} use the distance spectral radius to infer the extent of branching and model boiling points of an alkane.
 Therefore, maximizing or minimizing the distance spectral radius over a given class of graphs is of great interest and significance.
 Recently, the maximal or the minimal distance spectral radius of a given class of graphs has been studied extensively; see, e.g., 
 \cite{bose,ili,liu,nath,ste,yu2,yu1,zhang1,zhou1,zhou2}.

Recall that the {\it edge connectivity} of a connected graph is the minimum number of edges whose removal disconnects the graph.
For convenience, denote by $\mathcal{G}_n^r$ the set of all connected graphs of order $n$ and edge connectivity $r$.
Clearly, $1 \leq r \leq n-1$, and $\mathcal{G}_n^{n-1}$ consists of the unique graph $K_n$, where $K_n$ denotes a complete graph of order $n$.
Let $K(p,q)(p \geq q \geq 1)$ be a graph obtained from $K_p$ by adding a vertex together with edges joining this vertex to $q$ vertices of $K_p$.
Surely $K(n-1,r)\in \mathcal{G}_n^r$.
In this paper we prove that $K(n-1,r)$ is the unique graph with minimum distance spectral radius in $\mathcal{G}_n^r$, where $1 \le r \le n-2$.

\section{Main Results}
Given a graph $G$ on $n$ vertices, a vector $x \in \mathbb{R}^n$ is considered as a function defined on $G$,
  if there is a 1-1 map $\varphi$ from $V(G)$ to the entries of $x$; simply written $x_u=\varphi(u)$ for each $u \in V(G)$.
If $x$ is an eigenvector of $D(G)$, then it is naturally defined on $V(G)$, i.e. $x_u$ is the entry of $x$ corresponding to the vertex $u$.
One can find that
   $$x^{T}D(G)x =  \sum_{u,v \in V(G)}\dis_{uv} x_u x_v,\eqno(2.1)$$
  and $\la$ is an eigenvalue of $D(G)$ corresponding to the eigenvector $x$ if and only if $x \neq 0$ and
   $$\la x_v = \sum_{u \in V(G)} \dis_{vu} x_u,  \hbox{~ for each vertex ~} v \in V(G).  \eqno(2.2)$$
In addition, for an arbitrary unit vector $x \in {\mathbb R}^n$,
   $$x^TD(G)x\leq \rho(G), \eqno(2.3)$$
with the equality holds if and only if $x$ is an eigenvector of $D(G)$ corresponding to $\rho(G)$.

The following lemma is an immediate consequence of Perron-Frobenius theorem.

\begin{lemma} \label{edge}
Let $G$ be a connected graph with $u,v \in V(G)$. If $uv \notin E(G)$, then $\rho (G) > \rho (G+uv)$. If $uv \in E(G)$ and $G-uv$ is also connected, then $\rho (G) < \rho (G-uv)$.
 \end{lemma}

By Lemma \ref{edge}, for a connected graph $G$ on $n$ vertices, we have $\rho(G) \geq \rho(K_n)= n-1$, with equality holds if and only if $G=K_n$; and $\rho(G)\leq \rho(T_G)$, with equality holds if and only if $G=T_G$, where $T_G$ is a spanning tree of $G$.

Let $G$ be a graph and let $v$ be a vertex of $G$.
Denote by $N(v)$ the set of neighbors of $v$ in $G$, and by $d_v$ the degree of $v$ in $G$ (i.e. the cardinality of $N(v)$).

\begin{lemma} \label{inclu}
Let $G$ be a connected graph containing two vertices $u,v$, and let $x$ be a Perron vector of $D(G)$.\\
{\em(1)} If $N(u)\backslash \{v\}\subseteq N(v)\backslash \{u\}$, then $x_u \ge x_v$, with strict inequality if $N(u)\backslash \{v\}\subsetneq N(v)\backslash \{u\}$.\\
{\em(2)} If $N(u)\backslash \{v\} = N(v)\backslash \{u\}$, then $x_u = x_v$.
 \end{lemma}

{\it Proof:}
The second assertion follows from the first or can be found in \cite{liu}. So we only prove the assertion (1).
From (2.2), we have
$$\rho(G)x_u=\dis_{uv}x_v+\sum_{w \in V(G)\backslash \{u,v\}}\dis_{uw} x_w,\eqno (2.4)$$
$$ \rho(G)x_v=\dis_{vu}x_u+\sum_{w \in V(G)\backslash \{u,v\}}\dis_{vw} x_w. \eqno (2.5)$$
Since $N(u)\backslash \{v\}\subseteq N(v)\backslash \{u\}$, for each $w \in V(G)\backslash \{u,v\}$, we have
$$\dis_{uw} \geq \dis_{vw}, \eqno(2.6)$$
and hence
$$\sum \limits_{w \in V(G)\backslash \{u,v\}}\dis_{uw} x_w \ge \sum \limits_{w \in V(G)\backslash \{u,v\}}\dis_{vw} x_w. \eqno (2.7)$$
  By (2.4), (2.5) and (2.7), we get
$(\rho(G)+ \dis_{uv})x_u  \ge (\rho(G) +\dis_{uv})x_v$. So $x_u \ge x_v$.

If $N(u)\backslash \{v\}\subsetneq N(v)\backslash \{u\}$, then there exists a vertex $w \in (N(v)\backslash \{u\}) \setminus (N(u)\backslash \{v\})$ such that
$\dis_{uw}> \dis_{vw}=1$. So the inequality (2.6) is strict for some vertex $w$ and hence (2.7) holds strictly, which implies $x_u > x_v$.
\hfill$\blacksquare$

Let $G \in \mathcal{G}_n^r$. Then each vertex $v$ of $G$ holds $d_v \ge r$.
If there exists some vertex $v$ of $G$ with $d_v=r$, we have the following result immediately.

\begin{lemma} \label{degr}
Let $G \in \mathcal{G}_n^r\;(1 \le r \le n-2)$, which contains a vertex of degree $r$.
Then $\rho(G) \ge \rho(K(n-1,r))$, with equality if and only if $G=K(n-1,r)$.
\end{lemma}

{\it Proof}\ \ Let $v$ be a vertex of $G$ such that $d_v=r$.
Adding all possible edges within the subgraph of $G$ induced by the vertices of $V(G) \setminus \{v\}$,
we will arrive at a graph $G'$, which is isomorphic to $K(n-1,r)$.
If $G \ne G'$, then $\rho(G) > \rho(G')=\rho(K(n-1,r))$ by Lemma \ref{edge}.
The result follows.
\hfill$\blacksquare$

In the following we discuss the graph $G \in \mathcal{G}_n^r$ each vertex of which has degree greater than $r$.
We will formulate two lemmas about the behaviors of the distance spectral radius under some graph transformations,
and then establish the main result of this paper.

\begin{lemma} \label{per1}
Let $G$ be a graph obtained from $K_{n_1} \cup K_{n_2}$ by adding $r\;(\ge 1)$ edges between $u_1$ and $v_1,v_2,\ldots,v_r$, where $V(K_{n_1})=\{u_1,u_2,\ldots,u_{n_1}\}$, $V(K_{n_2})=\{v_1,v_2,\ldots,v_{n_2}\}$, $\min\{n_1,n_2\} \ge r+2$.
 Let $\tilde{G}$ be the graph obtained from $G$ by deleting the edges of $K_{n_1}$ incident to $u_1$ and adding all possible edges
 between the vertices of $V(K_{n_1})\backslash \{u_1\}$ and those of $V(K_{n_2})$.
 Then $\rho(G) > \rho (\tilde{G})$.
\end{lemma}

\begin{center}
\vspace{3mm}
\includegraphics[scale=.5]{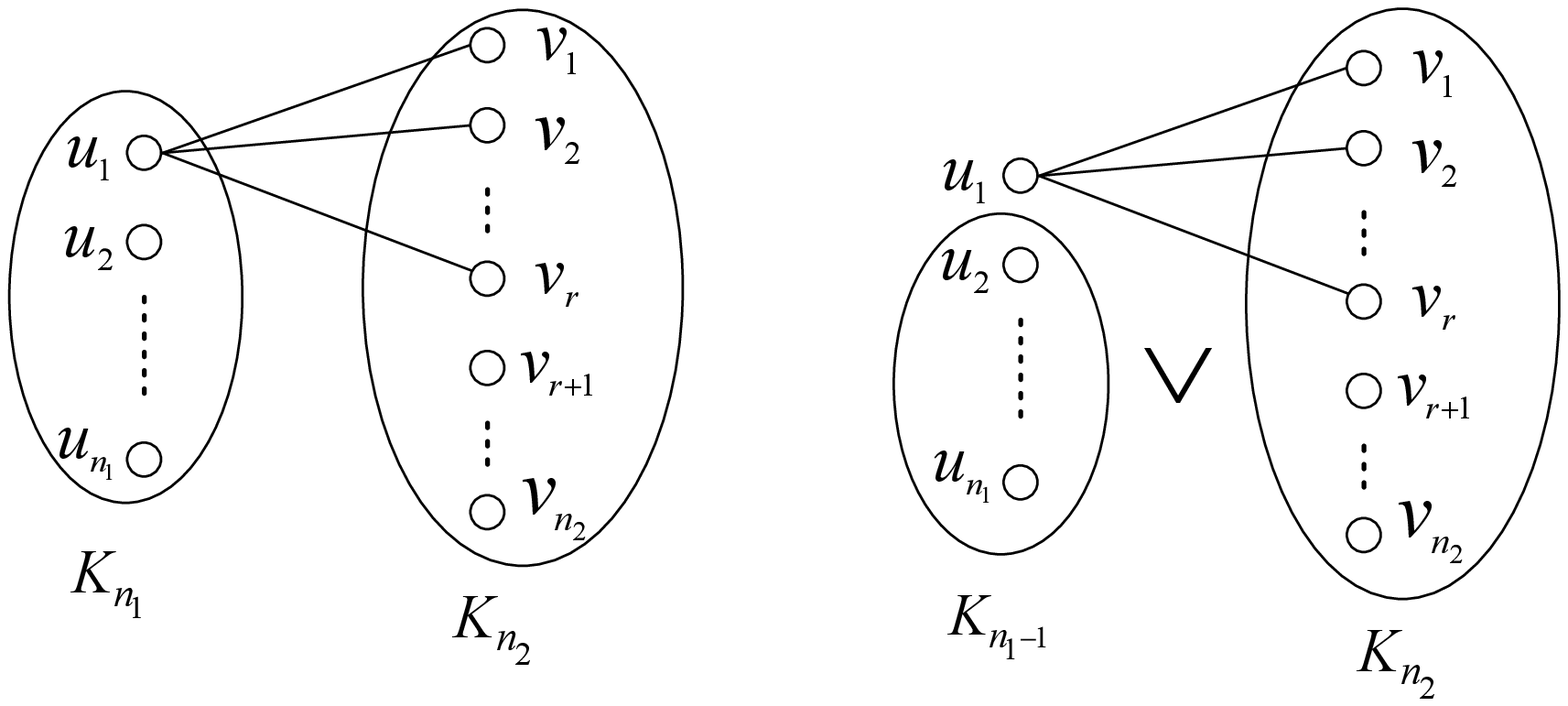}
\vspace{3mm}

{\small Fig. 2.1 \ \ The graphs $G$ (left) and $\tilde{G}$ (right) in Lemma \ref{per1}, \\
~~~~where $\vee$ means joining each vertex of $K_{n_1-1}$
and each of $K_{n_2}$}
\end{center}

{\it Proof:} \ \ Arrange in order the vertices of $\tilde{G}$ as $u_1,u_2,\ldots,u_{n_1},v_1,\ldots,v_r,v_{r+1},\ldots,v_{n_2}$.
Let $x$ be the unit Perron vector of $D(\tilde{G})$.
By Lemma \ref{inclu}, $x$ may be written as 
$$x=(x_1,\underbrace{x_2,\ldots,x_2}\limits_{n_1-1},\underbrace{x_3,\ldots,x_3}\limits_r,\underbrace{x_2,\ldots,x_2}\limits_{n_2-r})^T,\eqno (2.8)$$
where $x_3 < x_2 < x_1$.
Notice that the transformation from $G$ to $\tilde{G}$ leads to the distance between $u_1$ and $u_i \;(i=2,\ldots,n_1)$ increasing by 1,
the distance between $u_i \;(i=2,\ldots,n_1)$ and $v_j \;(j=1,\ldots,r)$ decreasing by 1,
and the distance between $u_i \;(i=2,\ldots,n_1)$ and $v_j \;(j=r+1,\ldots,n_2)$ decreasing by 2, while the distance between any other two vertices having no change.
Thus by (2.1) and (2.8),
$$\begin{array}{lcl} \displaystyle
x^TD(G)x-x^TD(\tilde{G})x& = &-2\sum\limits_{i=2,\ldots,n_1}x_{u_1}x_{u_i}+ 2 \sum\limits_{i=2,\ldots,n_1,\atop j=1,\ldots, r} x_{u_i}x_{v_j}+
  4 \sum\limits_{i=2,\ldots,n_1,\atop j=r+1,\ldots, n_2} x_{u_i}x_{v_j}\\
  &=& -2(n_1-1)x_1x_2+2r(n_1-1)x_2x_3+4(n_1-1)(n_2-r)x_2^2 \\
                   & = &2(n_1-1)x_2[-x_1+rx_3+2(n_2-r)x_2].
\end{array} \eqno (2.9)$$
Considering (2.2) on the vertex $u_1$ of $\tilde{G}$, we get
$$\rho(\tilde{G})x_1=\rho(\tilde{G})x_{u_1} = \sum_{i=1}^r x_{v_i}+\sum_{j=r+1}^{n_2}2x_{v_j}+\sum_{k=2}^{n_1}2x_{u_k}=rx_3+2(n_1+n_2-r-1)x_2.\eqno(2.10)$$
Noting that $\rho(\tilde{G})> n_1+n_2-1$, from (2.10) we have
$$x_1 = \frac{1}{\rho(\tilde{G})}[rx_3+2(n_1+n_2-r-1)x_2]<rx_3+2(n_2-r)x_2.\eqno(2.11)$$
By (2.9) and (2.11), we get $x^TD(G)x-x^TD(\tilde{G})x > 0.$
So according to (2.3) we get
$$\rho(G)\geq x^TD(G)x > x^TD(\tilde{G})x = \rho (\tilde{G}).$$
\hfill$\blacksquare$

\begin{lemma} \label{per2}
Let $G$ be a graph obtained from $K_{n_1} \cup K_{n_2}$ by adding $t \;(\ge 1)$ edges between $u_1$ and $v_1,v_2,\ldots,v_t$
and $r-t \;(\ge 1)$ edges between some vertices of $V(K_{n_1})\setminus \{u_1\}$ and some vertices of $V(K_{n_2})$,
where $V(K_{n_1})=\{u_1,u_2,\ldots,u_{n_1}\}$, $V(K_{n_2})=\{v_1,v_2,\ldots,u_{n_2}\}$, $\min\{n_1,n_2\} \ge r+2$.
 Let $\tilde{G}$ be a graph obtained from $G$ by deleting $n_1-(r+1-t)$ edges of $K_{n_1}$ between $u_1$ and $u_i$ for $i=2,\ldots,n_1-(r-t)$, and adding all possible edges between the vertices of $V(K_{n_1})\backslash \{u_1\}$ and those of $V(K_{n_2})$. Then $\rho(G) > \rho(\tilde{G})$.
\end{lemma}

\begin{center}
\vspace{3mm}
\includegraphics[scale=.5]{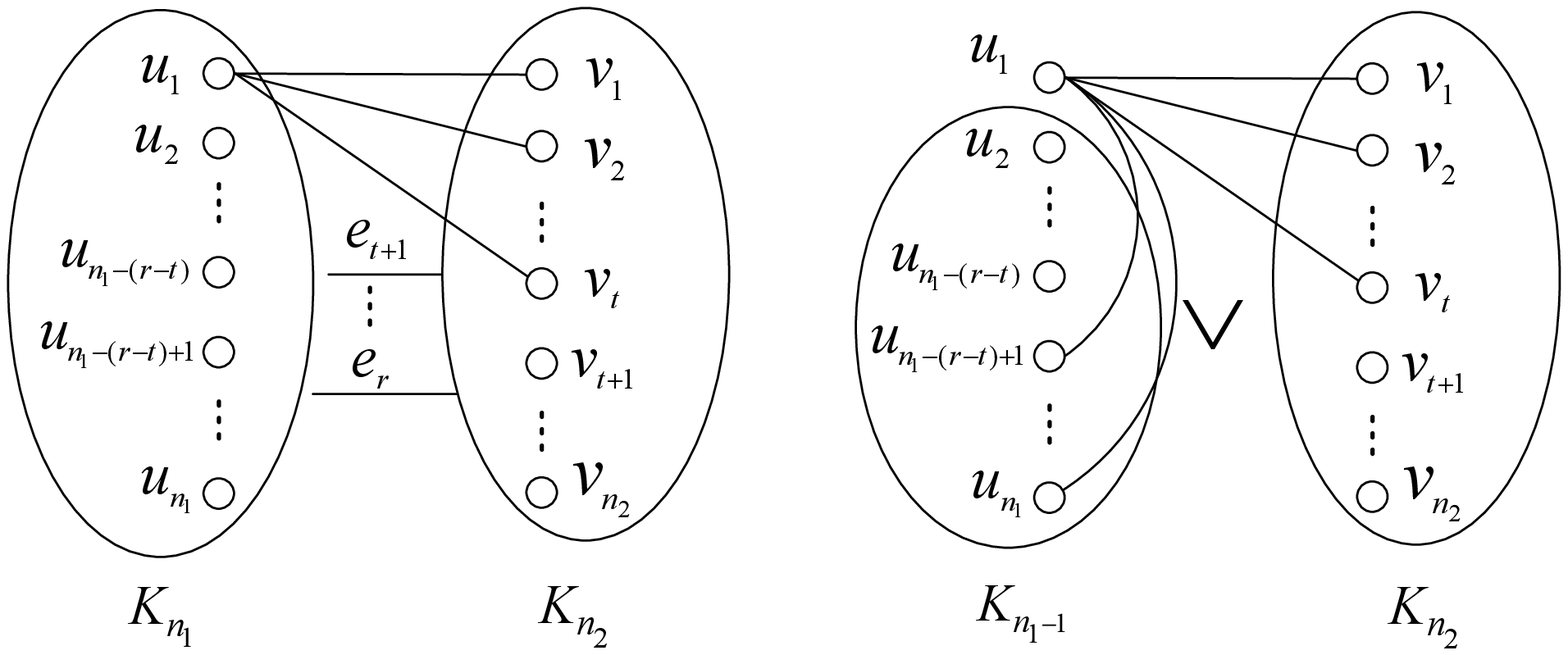}
\vspace{3mm}

{\small Fig. 2.2 \ \ The graphs $G$ (left) and $\tilde{G}$ (right) in Lemma \ref{per2}, \\
~~~~where $\vee$ means joining each vertex of $K_{n_1-1}$ and each of $K_{n_2}$}
\end{center}

{\it Proof:}\ \
 Arrange in order the vertices of  $V(\tilde{G})$ as $u_1,u_2,\ldots,u_{n_1-(r-t)}$, $u_{n_1-(r-t)+1},\ldots,u_{n_1}$, $v_1,\ldots,v_t$, $v_{t+1},\ldots,v_{n_2}$.
 Let $x$ be the unit Perron vector of $D(\tilde{G})$. By Lemma \ref{inclu}, $x$ may be written as $$x=(x_1,\underbrace{x_2,\ldots,x_2}\limits_{n_1-(r+1-t)},\underbrace{x_3,\ldots,x_3}\limits_{r-t},
 \underbrace{x_3,\ldots,x_3}\limits_t,\underbrace{x_2,\ldots,x_2}\limits_{n_2-t})^T,\eqno (2.12)$$
where $x_3 < x_2 < x_1$.
Let
$$ U_1=\{u_2,\ldots,u_{n_1-(r-t)}\}, U_2=\{u_{n_1-(r-t)+1},\ldots,u_{n_1}\},
V_1=\{v_1,\ldots,v_t\}, V_2=\{v_{t+1},\ldots,v_{n_2}\}.$$
Assume that in the graph $G$ there are $r_{ij}$ edges between $U_i$ and $V_j$ for $i,j=1,2$.
Surely, $r_{11}+r_{12}+r_{21}+r_{22}=r-t$.
Denote by $F$ the set of order pairs $(u,v)$ such that $uv$ is an edge of $G$, where $u \in U_1 \cup U_2, v \in V_1 \cup V_2$,
and by $(U_i,V_j)$ the set of order pairs $(u,v)$, where $u \in U_i, v \in V_j, i,j=1,2$.
Then by (2.1)
$$\begin{array}{lcl}
\frac{1}{2}[x^TD(\tilde{G})x-x^TD(G)x] & =& \sum\limits_{u \in U_1} x_{u_1}x_u -\sum\limits_{(u,v) \in (U_1,V_1) \setminus F}x_ux_v
  -\sum\limits_{(u,v) \in (U_1,V_2) \setminus F}\delta_{uv}x_ux_v \\
  & ~ &-\sum\limits_{(u,v) \in (U_2,V_1) \setminus F}x_ux_v
   -\sum\limits_{(u,v) \in (U_2,V_2) \setminus F}\delta_{uv}x_ux_v,
\end{array} \eqno (2.13)$$
where $\delta_{uv}=2$ if $u,v$ has distance $3$ in the graph $G$, and $\delta_{uv}=1$ otherwise.
By (2.12) and (2.13), and taking $\delta_{uv}=1$, we have
$$\begin{array}{lcl}
\frac{1}{2}[x^TD(\tilde{G})x-x^TD(G)x] & \le & [n_1-1-(r-t)]x_1x_2 -\{[n_1-1-(r-t)]t-r_{11}\}x_2x_3 \\
  & ~ &-\{[n_1-1-(r-t)](n_2-t)-r_{12}\}x_2^2 -[(r-t)t-r_{21}]x_3^2\\
  & ~ &-[(r-t)(n_2-t)-r_{22}]x_2x_3\\
  &=& [n_1-1-(r-t)]x_1x_2 -\{[n_1-1-(r-t)]t+(r-t)(n_2-t)\}x_2x_3 \\
  &~& -\{[n_1-1-(r-t)](n_2-t)\}x_2^2 -(r-t)tx_3^2 \\
  &~& + (r_{11}+r_{22})x_2x_3+r_{12}x_2^2+r_{21}x_3^2\\
  & \le & [n_1-1-(r-t)]x_1x_2 -\{[n_1-1-(r-t)]t+(r-t)(n_2-t)\}x_2x_3 \\
  &~& -\{[n_1-1-(r-t)](n_2-t)\}x_2^2 -(r-t)tx_3^2 + (r-t)x_2^2\\
  &=& [n_1-1-(r-t)]x_1x_2 -\{[n_1-1-(r-t)]t+(r-t)(n_2-t)\}x_2x_3 \\
  &~& -\{[n_1-1-(r-t)](n_2-t)-(r-t)\}x_2^2 -(r-t)tx_3^2
\end{array} \eqno (2.14)$$
Considering (2.2) on the vertex $u_1$ of $\tilde{G}$, we get
$$\rho(\tilde{G})x_1 = rx_3+2(n_1+n_2-r-1)x_2.$$
So $x_1 = \frac{1}{\rho(\tilde{G})}[rx_3+2(n_1+n_2-r-1)x_2]<\frac{1}{n_1+n_2-1}[rx_3+2(n_1+n_2-r-1)x_2]$.
Therefore,
$$\begin{array}{lcl}
\frac{1}{2}x^T[D(\tilde{G})-D(G)]x &<&\frac{n_1-1-(r-t)}{n_1+n_2-1}rx_2x_3+\frac{[n_1-1-(r-t)]\cdot2(n_1+n_2-r-1)}{n_1+n_2-1}x_2^2\\
                                       &~&-\{[n_1-1-(r-t)]t+(r-t)(n_2-t)\}x_2x_3\\
                                       &~& -\{[n_1-1-(r-t)](n_2-t)-(r-t)\}x_2^2 -(r-t)tx_3^2
\end{array} \eqno (2.15)$$
Let $a,b,c$ be the coefficients of $x_2^2,x_2x_3,x_3^2$ in (2.15), respectively.
Noting  that $\min\{n_1,n_2\} \ge r+2$, we have
$$\begin{array}{lcl}
a&=&2[n_1-1-(r-t)](1-\frac{r}{n_1+n_2-1})-\{[n_1-1-(r-t)](n_2-t)-(r-t)\}\\
 &=&2[n_1-1-(r-t)](1-\frac{r}{n_1+n_2-1}-\frac{n_2-t}{2})+(r-t)\\
 &<&2[n_1-1-(r-t)](1-\frac{n_2-t}{2})+(r-t)\\
 &\leq&2[n_1-1-(r-t)]\frac{2-n_2+t}{2}+(n_2-2-t)\\
 &=&-(n_2 -t-2)[n_1-2-( r-t)] < 0\\
b&=&[n_1-1-(r-t)](\frac{r}{n_1+n_2-1}-t)-(r-t)(n_2-t) < 0,\\
c&=&-(r-t)t < 0\\
\end{array} $$
Thus $x^TD(\tilde{G})x-x^TD(G)x<0$, and hence $\rho(G)\geq x^TD(G)x > x^TD(\tilde{G})x = \rho (\tilde{G}).$
\hfill$\blacksquare$

\begin{lemma} \label{order}
Let $G$ be a connected graph, and let $E_c$ be an edge cut set of $G$ of size $r \;(\ge 1)$ such that $G-E_c = K_{n_1} \cup K_{n_2}$, where $n_1 + n_2 = n$.
If $d_v > r$ for each vertex $v \in V(G)$, then $n_1 \ge r+2, n_2 \ge r+2$.
\end{lemma}

{\it Proof:}\ \ 
 If $n_1 \le r$, then there exists a vertex $u$ of $K_{n_1}$ such that
 $$d(u) \leq n_1-1+ \frac{r}{n_1} \le (n_1-1)\frac{r}{n_1} + \frac{r}{n_1} = r,$$ a contradiction.
If $n_1 = r+1$, then there exists a vertex $w$ not incident with any edges of $E_c$, which implies $d(u) = r$, also a contradiction.
The discussion for the assertion on $n_2$ is similar. \hfill$\blacksquare$

\begin{theorem}
For each $r=1,2,\ldots,n-2$, the graph $K(n-1,r)$ is the unique graph with minimum distance spectral radius in $\mathcal{G}_n^r$.
\end{theorem}
{\it Proof:} Let $G$ be a graph that attains the minimum distance spectral radius in $\mathcal{G}_n^r$.
Note that each vertex of $G$ has degree not less than $r$.
If there exists a vertex $u$ of $G$ with degree $r$, by Lemma \ref{degr}, $\rho(G) \ge \rho(K(n-1,r))$, with equality if and only if $G=K(n-1,r)$.
So the result follows in this case.

Next we assume all vertices of $G$ have degrees greater than $r$.
Let $E_c$ be an edge cut set of $G$ containing $r$ edges, and let $G_1,G_2$ be two components of $G-E_c$ with order $n_1,n_2$ respectively.
We assert $G_1=K_{n_1}$ and $G_2=K_{n_2}$; otherwise adding all possible edges within $G_1,G_2$ we would get a graph with smaller distance spectral radius by Lemma \ref{edge}.
By Lemma \ref{order}, $n_1 \ge r+2, n_2 \ge r+2$.
Let $u_1$ be a vertex of $G_1$ such that $u_1$ joins $t$ vertices of $G_2$, where $1 \le t \le r$.
If $t=r$, by Lemma 3.2 there exists a graph $\tilde{G} \cong K(n-1,r)$ such that $\rho(G) > \rho (\tilde{G})$.
If $1 \le t<r$, by Lemma 3.3 there also exists a graph $\tilde{G} \cong K(n-1,r)$, such that $\rho(G) > \rho (\tilde{G})$.
This completes the proof.\hfill$\blacksquare$


\small
\begin{thebibliography}{90}
\bibitem{bal} A. T. Balaban, D. Ciubotariu, M. Medeleanu, Topological indices and real number vertex invariants based on graph eigenvalues or eigenvectors, {\it J. Chem. Inf. Comput. Sci.}, 31(1991), 517-523.
\bibitem{bose} S. S. Bose, M. Nath, S. Paul, Distance spectral radius of graphs with $r$ pendent vertices, {\it Linear Algebra Appl.}, 435(2011), 2826-2836.

\bibitem{ede} M. Edelberg, M. R. Garey and R. L. Graham, On the distance matrix of a tree, {\it Discrete Math.}, 14(1976), 23-29.

\bibitem{gra3} R. L. Graham and L. Lovasz, Distance matrix polynomials of trees, {\it Adv. Math.}, 29(1978), 60-88.

\bibitem{gra1} R. L. Graham and H. O. Pollak, On the addressing problem for loop switching, {\it Bell Sys. Tech. J.}, 50(1971), 2495-2519.

\bibitem{gra2} R. L. Graham and H. O. Pollak, On embedding graphs in squashed cubes, In: {\it Graph Theory and Applications}, Springer, Berlin, 1973, 99-110.

\bibitem{gut}I. Gutman, M. Medeleanu, On structure-dependence of the largest eigenvalue of the distance matrix of an alkane, {\it Indian J. Chem. A}, 37(1998), 569-573.

\bibitem{hos} H. Hosoya, M. Murakami and M. Gotoh, Distance polynomial and characterization of a graph, {\it Natur. Sci. Rep. Ochanomizu Univ.}, 24(1973), 27-34.

\bibitem{ili} A. Ili\'c, Distance spectral radius of trees with given matching number, {\it Discrete Appl. Math.}, 158(2010), 1799-1806.

\bibitem{liu} Z. Z. Liu, On the spectral radius of the distance matrix, {\it Appl. Anal. Discrete Math.}, 4(2)(2010), 269-277.

\bibitem{nath} M. Nath, S. Paul, On the distance spectral radius of bipartite graphs, {\it Linear Algebra Appl.}, 436(2012), 1285-1296.

\bibitem{rou} D. H. Rouvray, The search for useful topological indices in chemistry, {\it Amer. Scientist}, 61(1973), 729-735.

\bibitem{ste} D. Stevanovi\'c, A. Ili\'c, Distance spectral radius of trees with fixed maximum degree, {\it Electron. J. Linear Algebra}, 20(2010), 168-179.

\bibitem{yu2} G. L. Yu, H. C. Jia, H. L. Zhang, J. L. Shu, Some graft transformations and its application on a distance spectrum, {\it Appl. Math. Letters}, 25(2012), 315-319.

\bibitem{yu1} G. L. Yu, Y. R. Wu, J. L. Shu, Some graft transformations and its application on a distance spectrum, {\it Discrete Math.}, 311(2011), 2117-2123.


\bibitem{zhang1} X. L. Zhang, C. Godsil, Connectivity and minimal distance spectral radius, {\it Linear Multilinear Algebra}, 59(2011), 745-754.


\bibitem{zhou1} B. Zhou, On the largest eigenvalue of the distance matrix of a tree, {\it MATCH Commun. Math. Comput. Chem.}, 58(2007), 657-662.

\bibitem{zhou2} B. Zhou, N. Trinajisti\'c, On the largest eigenvalue of the distance matrix of a connected graph, {\it Chem. Phys. Lett.}, 447(2007), 384-387.




\end{thebibliography}
\end{document}